\input amstex
\documentstyle{amsppt}
\NoBlackBoxes
\magnification=1200
\hcorrection{.25in}
\advance\vsize-.75in
\topmatter
\title\nofrills OBSTRUCTIONS TO SHELLABILITY
\endtitle
\author Michelle L. Wachs
\endauthor
\address Department of Mathematics, University of Miami, Coral Gables,
FL 33124
\endaddress
\email wachs\@math.miami.edu \endemail
\thanks Research supported in part by  NSF
grant DMS 9311805 at the University of Miami and by NSF grant DMS-9022140
at MSRI.
\endthanks
\abstract We consider a simplicial complex generaliztion of a result
of Billera and Meyers   that every nonshellable poset contains the
smallest nonshellable poset as an induced subposet. We prove that every
nonshellable $2$-dimensional simplicial complex contains a nonshellable
induced subcomplex with less than $8$ vertices.   We also establish
CL-shellability of interval orders and as a consequence obtain a formula for
the Betti numbers of any interval order.

\endabstract
\endtopmatter
\document
\def \pure{\operatorname{pure}}
\def \lk{\operatorname{lk}}

A recent result of Billera and Meyers \cite{BM} implies 
 that every
nonshellable poset contains as an induced subposet the 4 element poset
$Q$  consisting of two disjoint 2 element chains.   (Throughout this
paper shellability refers to the general notion of nonpure shellability
 introduced in
\cite {BW2}.)
  Note that $Q$ is the  nonshellable poset with the fewest number of
elements. 
 Of course a
shellable poset can also contain $Q$; eg.,  the lattice of
subsets of a 3 element set.  So the condition of not containing $Q$ as
an induced subposet  is only sufficient for shellability; it does not
characterize shellability.  It is however a well-known
characterization of  a class of posets called interval orders and the
question of whether all interval  orders are shellable is what
Billera and Meyers were considering in the first
place.  

 In this note we
suggest a way to generalize  the poset result to general
simplicial  complexes.   We also give a simple proof of the
poset result and   prove the stronger result that any
poset that does not contain
$Q$ as an induced subposet is
$CL$-shellable.  This yields a recursive formula for the Betti
numbers of the poset.  

We assume familiarity with the general theory of shellability
\cite{BW2} \cite{BW3}.  All notation and terminology used here is defined in
\cite{BW2 }and \cite{BW3}.

The most simple minded conjecture one could make is
that every nonshellable simplicial complex contains the induced
subcomplex  consisting of edges
$\{a,b\}$ and $\{c,d\}$, where $a,b,c,d$ are distinct vertices. 
  A simple counterexample
 is given by the 5 vertex simplicial complex
consisting of facets
$\{a,b,c\},
\{c,d,e\}, \{a,d\}$.  Indeed the situation for simplicial complexes
turns out to be much more complicated than it is for posets.

The 
 most natural thing to do next is to look for other
``obstructions'' to simplicial complex shellability.  Is there a finite
list?   Below we see that the answer is no. Define an {\it obstruction}
to shellability to be a nonshellable simplicial complex  all of whose
proper  induced  subcomplexes are shellable.  The 4 and 5 vertex
simplicial complexes given above are examples of one and
two dimensional obstructions, respectively.  The following observation
was made by Stanley \cite{S}.

\proclaim{Proposition 1} For every positive integer $d$ there is an
obstruction to shellability of dimension $d$.
\endproclaim

\demo{Proof}  Let $K$ be the
$(d-1)$-skeleton  of the simplex on vertex set $[d+3]=
\{1,2,\dots,d+3\}$ together with  two
$d$-dimensional faces
$\{1,2,\dots,d+1\}$ and $
\{3,4,\dots,d+3\}$.  We claim that $K$ is a $d$-dimensional
obstruction.  If
$K$ were shellable then by the Rearrangement Lemma of \cite {BW2}
there  would be a shelling order in which the maximal dimensional
facets come first; namely 
$\{1,2,\dots,d+1\}$ and
$\{3,4,\dots,d+3\}$ come first.  But this is impossible because these
two facets intersect in a face of dimension $d-2$.  Hence $K$ is not
shellable. 

 Every proper induced subcomplex of $K$ is either a simplex or
consists of a single $d$-face in a $(d+1)$-simplex together with the
$(d-1)$ skeleton of  the $(d+1)$-simplex.  Certainly the simplex is
shellable. Let $J$ be the $(d-1)$-skeleton of the simplex on vertex
set $ [d+2]$ together with the face
$\{1,2,\dots,d+1\}$.  It is easy to see that lexicographical order on
the facets of $J$ is a shelling of $J$. (Lexicographical order on
subsets of $[d+2]$ is defined by $\{a_1 \le \dots \le a_k\} \le \{b_1
\le
\dots
\le b_j
\}$ if the word $a_1 \cdots a_k$ is less than the word $b_1 \cdots
b_j$ in lexicographical order.) Therefore all proper induced
subcomplexes of $K$ are shellable.\hfill\qed
\enddemo

We now consider the following problem.

\proclaim{Problem} Determine whether or not there is a finite number of
$d$-dimensional obstructions to shellability for each $d$.  If so, find
bounds on the number of vertices that a $d$-dimensional obstruction can have.
\endproclaim 


In this paper we solve this problem only for dimensions $d =1,2$ and we
leave open the problem for general $d$.

\proclaim{Proposition 2}  The only 1-dimensional obstruction to
shellability is the complex $J$ generated by facets $\{a,b\}$,
$\{c,d\}$ where
$a,b,c,d$ are distinct.  
\endproclaim

\demo{Proof}  A $1$-dimensional simplicial complex is shellable if
and only if it has at most one connected component with more than $1$
vertex.  If $K$ is a nonshellable $1$-dimensional
simplicial complex then let $\{a,b\}$ be an edge in one
 component of $K$ and let $\{c,d\}$ be an edge in
another  component.  The subcomplex induced by
$a,b,c,d$ is $J$. Hence $K$ is an obstruction if and only if
$K=J$. \hfill\qed
\enddemo

Already in dimension 2 the situation is much more complicated.
We use the following notation:  For any subset $U$ of $V$ and
simplicial complex
$K$ on vertex set $V$, let $K(U)$ be the subcomplex of
$K$ induced by $U$. Also let the pure part of $K$, denoted
$\pure(K)$, be the subcomplex of $K$ generated by the facets of
maximum dimension.  For $v \in V$, the link of $v$ in $K$ is denoted
$\lk_K(v)$ and is defined to be $\{F \in K \mid F \cup \{v\} \in K
\text{ and } v \notin F\}$.  For any $v \in
V$ and subcomplex $J$ of $\lk_K(v)$,  the join of $v$ and $J$ is
denoted
$v*J$ and is defined to be
$\{F \in K \mid v \in F \text{ and } F \setminus \{v\} \in J\}$.
The ith reduced simplicial homology of $K$ over the ring of integers is
denoted by $\tilde H_i(K)$.

\proclaim{Theorem 3} There are no 2-dimensional obstructions with more
than 7 vertices.
\endproclaim

\demo{Proof}
Let $K$ be a 2-dimensional simplicial complex on vertex set $V$ where
$|V| > 7$.  Assume all induced proper subcomplexes are shellable.  We shall
show that
$K$ is shellable  by showing that $\pure(K)$ is
shellable and the 1-skeleton of $K$ is connected (except for isolated
points).  That the 1-skeleton is connected follows immediately from the
fact that no induced subcomplex consists only of   a pair of
disjoint edges.  

To prove that $\pure(K)$ is shellable, choose any vertex $v$ of
$\pure(K)$.   Let
$$K_1 = \pure(K(V\setminus\{v\})).$$  Since the pure part of a shellable
complex is shellable (by the Rearrangement Lemma of \cite{BW2}),
$K_1$ is shellable.  Let $$K_2 =v*\pure(\lk_{K}(v)).$$  We claim that
$\pure(\lk_{K}(v))$ is  a connected 1-dimensional complex.   If not there
would be distinct  vertices $a,b,c,d \in V \setminus
\{v\}$ such that edges
$\{a,b\}$ and
$\{c,d\}$ are in different components of $\pure(\lk_K(v))$. Since
$|V| > 5$, the induced subcomplex $K(\{v,a,b,c,d\})$ would be  shellable
which would imply that $\lk_{K(\{v,a,b,c,d\})}(v)$ is shellable
since the link of any vertex in a shellable complex is shellable
\cite{BW3}.  But this is impossible since
$\lk_{K(\{v,a,b,c,d\})}(v)$ has only two facets $\{a,b\}$ and $\{c,d\}$.
It follows from this claim that $K_2$ is shellable and $2$-dimensional.

Next we dispose of the special cases that  $K_1$ 
is 0 or 1-dimensional.   Clearly  $K_1$ can't be $0$-dimensional since $v$
belongs to a $2$-face. If
$K_1$ is 1-dimensional then
$\pure(K) = K_2$  which is shellable.  

Now we can assume that $K_1$ and $K_2$  are both shellable and
2-dimensional. Note that then
$$\pure(K) = K_1 \cup K_2.$$    Let $$A =
K_1
\cap K_2.$$

  We shall show that $A$ is connected and 1-dimensional.
Suppose it isn't.  Then either (1) $A = \{\emptyset\}$, (2) $A$ contains an
isolated point
 or (3) there are edges  in different components of $A$.  For the first
case, choose any
$a,,b,c,d,e$ such that
$\{a,b\} \in \pure(\lk_K(v))$ and $\{c,d,e\} \in K_1$. Since $|V| > 6$,
$K(\{v,a,b,c,d,e\})$ is shellable.  It follows that $\{v,a,b\}$ and
$\{c,d,e\}$ cannot be the only 2-faces of $K(\{v,a,b,c,d,e\})$. Hence there
is a third 2-face $F$.  Note $v \notin F$ because otherwise one of the other
vertices of $F$ would be in $A$.  So $F \in K_1$.  Since either $a$ or $b$
is a vertex of $F$ as well as of $\pure(\lk_K(v))$, $a$ or $b$ is a vertex of
$K_1 \cap K_2$.  Hence $A$ cannot be $\{\emptyset\}$.  

For the second  case, let $a$ be the isolated point.  Then
$K_1$ and
$K_2$ contain 2-faces 
$\{a,c,d\}$ and $\{v,a,b\}$, respectively, which intersect only at $a$. 
Since
$|V| > 5$, $K(\{v,a,b,c,d\})$ is shellable.  This means
that there is a third 2-face in the induced subcomplex  that intersects
each of the 2-faces along  edges that  contain $a$.  If the third
2-face  contains $v$ then it is either $\{v,a,c\}$ or $\{v,a,d\}$. 
This implies that either $\{a,c\}$ or $\{a,d\}$ is in $A$, which
contradicts the fact that $\{a\}$ is a facet of $A$. Hence the third
2-face must be $\{a,b,c\}$ or $\{a,b,d\}$. It follows
that 
$\{a,b\}$ is a facet of $A$, which is still a contradiction.

For  the third case, suppose that $\{a,b\}$ and $\{c,d\}$ are edges in
different  components of $A$.  Let
$x,y
\in V
\setminus \{v\}$ be such that
$\{a,b,x\}$ and $\{c,d,y\}$ are facets of $K_1$.  Since  $|V| > 7$,
$J=K(\{v,a,b,c,d,x,y\})$ is shellable.  Let 
$$B=(v*\pure(\lk_J(v))) \cap
\pure(K(\{a,b,c,d,x,y\})).$$ 
Since $B$ is a subcomplex of $A$, $\{a,b\}$
and $\{c,d\}$ are in different components of $B$.  It follows that
$\tilde H_0(B) \ne 0$.  Since $v*\pure(\lk_J(v))$ is
contractible, we also have $\tilde H_i(v*\pure(\lk_J(v))) = 0$
for all $i$.  By the Rearrangement Lemma of
\cite{BW2}, $\pure(K(\{a,b,c,d,x,y\})$  is shellable since the
induced subcomplex $K(\{a,b,c,d,x,y\})$ is. Hence $\tilde
H_i(\pure(K(\{a,b,c,d,x,y\}))) = 0$ for $i \le 1$.  Since 
$$\pure(J) = (v*\pure(\lk_J(v))) \cup
\pure(K(\{a,b,c,d,x,y\})),$$
 by  Mayer-Vietoris 
we have that $\tilde H_1(\pure(J)) =\tilde H_0(B) \ne 0$.  This
contradicts the fact that $J$ is shellable and 2-dimensional.  Hence we
may conclude that 
$A$ is connected and 1-dimensional.

Since
$A$ and 
$\pure(\lk_K(v))$ are connected we can get a shelling of  
$\pure(\lk_K(v))$ by first listing the edges of $A$  and then listing the
remaining edges of
$\pure(\lk_K(v))$ so that each edge is connected to the previous
ones.  We claim that we can obtain a shelling of
$\pure(K)$ by first listing the facets of $K_1$  in the order given
by the shelling of
$K_1$  and then listing the facets of $K_2=v*\pure(\lk_K(v))$ in the order
indicated by the shelling of $\pure(\lk_K(v))$ in which the edges of $A$
come first.  Let
$F_1,F_2,\dots,F_n$ be the resulting ordered list of facets of
$\pure(K)$. 

 If
$F_i \in K_1$ then it is clear that  $(\bar F_1 \cup \bar F_2
\cup
\dots \cup \bar F_{i-1}) \cap \bar F_i$ is   pure
1-dimensional.  If
$F_i =
\{v,a,b\}$ where $\{a,b\} \in A$ then it is easy to see that  
$\{a,b\}$ is a facet of $(\bar F_1 \cup \bar F_2
\cup
\dots \cup \bar F_{i-1}) \cap \bar F_i$ and that for all but
the first $F_i$ in $v*A$, $\{v,a\}$ or $\{v,b\}$ is also a facet 
of $(\bar F_1 \cup \bar F_2
\cup
\dots \cup \bar F_{i-1}) \cap \bar F_i$. Hence $(\bar F_1 \cup
\bar F_2
\cup
\dots \cup \bar F_{i-1}) \cap \bar F_i$ is pure 1-dimensional
for all $F_i \in v*A$.  

Now suppose $F_i = \{v,a,b\}$ where
$\{a,b\}
\notin A$.  Clearly either $\{v,a\}$ or $\{v,b\}$ is  a facet of $
(\bar F_1
\cup
\bar F_2
\cup
\dots \cup \bar F_{i-1}) \cap \bar F_i$; assume without loss
of generality that $\{v,a\}$ is the facet.  If
$\{b\}$ is a facet of $ (\bar F_1 \cup
\bar F_2
\cup
\dots \cup \bar F_{i-1}) \cap \bar F_i$, then $b$ is in some
2-face $\{b,c,d\}$ of $K_1$.  Consider the induced subcomplex
$L=K(\{v,a,b,c,d\})$. Since
$\pure(L)$ is shellable,  $L$ contains one of the following 2-faces:
$\{a,b,c\}$, $\{a,b,d\}$, $\{v,b,c\}$ or $\{v,b,d\}$. The first two are
impossible since  $\{a,b\} \notin A$.  The last two facets
have to precede $\{v,a,b\}$ in the shelling because $\{b,c\}$ and
$\{b,d\}$ are in $A$. It follows that  $\{v,b\} \in (\bar F_1 \cup
\bar F_2 
\cup
\dots \cup \bar F_{i-1}) \cap \bar F_i$ which contradicts the
assumption that $\{b\}$ is a facet.  Hence $\{b\}$ is not a facet of $
(\bar F_1 \cup
\bar F_2
\cup
\dots \cup \bar F_{i-1}) \cap \bar F_i $ which implies that 
$ (\bar F_1 \cup
\bar F_2
\cup
\dots \cup \bar F_{i-1}) \cap \bar F_i$ is pure
1-dimensional. Therefore $F_1,F_2,\dots,F_n$ is indeed a shelling of
$\pure(K)$.
\hfill\qed
\enddemo

\proclaim{Corollary 4} Every nonshellable $2$-dimensional simplicial
complex  has a nonshellable induced subcomplex with $n$ vertices
where $4 \le n \le 7$.\hfill \qed
\endproclaim

\proclaim{Theorem 5} For each $n=5,6,7$,  there is a $2$-dimensional
obstruction with $n$ vertices. 
\endproclaim

\demo{Proof}  Let $M_n$ be the simplicial  complex on vertex set
$\{1,\dots ,n\}$ with facets $\{1,2,3\}$, $\{2,3,4\}$, $\dots$,
$\{n-2,n-1,n\}$, $\{n-1,n,1\}$, $ \{n,1,2\}$.  For  $n \ge 5$,
$M_n$ triangulates a cylinder when $n$ is even and a M\"obius strip
when $n$ is odd.  Hence $M_n$ is not shellable.  We leave it to the
reader to check that every induced proper subcomplex of $M_n$ is
shellable when $n \le 7$.\hfill\qed
\enddemo

The obstructions given in the proof of Theorem 5 are pseudomanifolds
with boundary.  We show next that a pseudomanifold without boundary
cannot be an obstruction in any dimension.

\proclaim{Lemma 6} Let $K$ be a simplicial complex on vertex
set $V$.  Suppose that $v \in V$ is such that $K(V \setminus \{v\})$
and $\lk_K(v)$ are shellable and no facet of $\lk_K(v)$ is a
 facet of $K(V \setminus \{v\})$. Then $K$ is shellable.
\endproclaim

\demo{Proof}  First list the facets of $K(V \setminus \{v\})$ in any
shelling order and then list the facets of $v* \lk_K(v)$ in the order
indicated by the shelling of $\lk_K(v)$. It is easy to see that this is
a shelling of $K$.\hfill\qed
\enddemo

\proclaim{Theorem 7}  Let $K$ be a  simplicial complex for which
every nonfacet face is contained in at least two facets.  Then
$K$ is not an obstruction.  Consequently,  there are no obstructions
that are pseudomanifolds (without boundary) or triangulations of
manifolds (without boundary). 
\endproclaim

\demo{Proof} The proof is by induction on $\dim K$.  When $\dim K =
1$ the result follows immediately from Proposition 2.  Suppose
$\dim K >1$ and every proper induced subcomplex of
$K$ is shellable. We will show that $K$ must also be shellable. Let
$V$ be the vertex set of $K$ and choose any $v \in V$.  Then
$K(V\setminus
\{v\})$ is  shellable. It follows from the fact that every
nonfacet face is in at least two facets of $K$ that no facet
of $\lk_K(v)$ is a facet of  $K(V\setminus
\{v\})$.  

To apply Lemma 6 we need only show that $\lk_K(v)$ is
shellable.  Note that the property that every nonfacet face is
contained in at least two facets, is inherited by  $\lk_K(v)$.  Hence
 if 
$\lk_K(v)$ is not shellable then by induction it  contains an
obstruction $\lk_K(v)(U)$, where $U \subsetneq V \setminus\{v\}$. We
have that 
$K(U\cup\{v\})$ is shellable since it is a  proper induced subcomplex
of $K$.  Since
$ \lk_K(v)(U)=\lk_{K(U\cup\{v\})}(v)
$ and any link in a shellable complex is shellable, we have that $
\lk_K(v)(U)$ is also shellable, contradicting the fact that $
\lk_K(v)(U)$ is an obstruction.  Therefore $\lk_K(v)$ is shellable
and by Lemma 6, $K$ is shellable. So $K$ is not an
obstruction.\hfill\qed
\enddemo

A ``pure'' version of Lemma~6 is used implicitly in Provan and Billera's
proof of shellability of matroid complexes \cite{PB}.  A {\it matroid
complex} is a simplicial complex for which all   induced subcomplexes are
pure.  Lemma~6 can, in fact, be used to prove the following stronger result.

\proclaim{Proposition 8} If every proper induced subcomplex of a simplicial
complex $K$ is pure then $K$ is shellable.
\endproclaim
\demo{Proof} The proof is by induction on the size of the vertex set $V$. 
Suppose that $K$ is not a simplex.
  Let
$F$ be any
$d$-dimensional facet of $K$ where $d=\dim K$. Choose $v \in V \setminus
F$.   Clearly $K(V \setminus \{v\})$ is pure
$d$-dimensional since it contains $F$.  It follows that no facet of 
$\lk_K(v)$ is a facet of  
$K(V
\setminus \{v\})$.  

By induction,
$K(V
\setminus
\{v\})$ is  shellable.  To apply Lemma~6, we need only show that  $\lk_K(v)$
is also shellable.   For  any $U \subsetneq V\setminus \{v\}$, $K(U \cup
\{v\})$ is pure.  Since 
$\lk_K(v)(U) =
\lk_{K(U
\cup \{v\})}(v)$ and any link in a pure complex is pure, we have that
$\lk_K(v)(U)$ is pure.   Hence every proper induced subcomplex of
$\lk_K(v)$ is pure. It follows by induction that $\lk_K(v)$ is
shellable.\hfil\qed
\enddemo

\demo {Remark}  Provan and Billera \cite{PB}  prove that matroid
complexes are  shellable by showing that they are vertex decomposible.  The
proof of Proposition 8 given here is a slight modification of the
Provan-Billera proof and also yields the  conclusion that $K$ is
vertex decomposible, but in the nonpure sense described in \cite{BW 3}.
\enddemo

Define an {\it obstruction to purity} to be a nonpure simplicial complex
for which all proper induced subcomplexes are pure.  Proposition~8
extends the Provan-Billera result from matroid complexes to  obstructions to
purity.   It turns out that there are really very few obstructions to purity
in each dimension.

\proclaim{Proposition 9} For each $d \ge 1$, every $d$-dimensional
obstruction to purity has exactly $d+2$ vertices.  Moreover there exists a
$d$-dimensional obstruction to purity for each $d$.
\endproclaim 

\demo{Proof} Let $K$ be a $d$-dimensional simplicial complex with vertex
set $V$.  Suppose  $|V| > d+2$ and all proper induced subcomplexes of $K$
are pure.  We will show that $K$ is also pure.  

Let $v \in V\setminus F$, where $F$ is any $d$-dimensional face of $K$.  It
follows that $K(V\setminus\{v\})$ is pure 
$d$-dimensional.  To show that
$K$ is pure we need only show that $\lk_K(v)$ is pure $(d-1)$-dimensional. 
Let $G$ be a face of $\lk_K(v)$.  Then since $K(V\setminus\{v\})$ is pure 
$d$-dimensional, $G$ is contained in some $d$-dimensional facet $H$ of
$K(V\setminus\{v\})$.  Since $H \cup \{v\} \ne V$, it follows that  $K(H \cup
\{v\})$ is pure
$d$-dimensional.  Since $G \cup \{v\} \in K(H \cup
\{v\})$, it follows that $G
\cup \{v\}$ is contained in a facet of dimension $d$ in $K(H \cup
\{v\})$.  This implies that
$G$ is contained in a $(d-1)$-dimensional face of $\lk_K(v)$, which
means that $\lk_K(v)$ is pure and $(d-1)$-dimensional.

An example of a $d$-dimensional obstruction to purity is  given by the
$(d-1)$-skeleton of the simplex on vertex set
$[d+2]$ together with the
$d$-dimensional face $[d+1]$.
\hfill\qed
\enddemo

Lemma 6 also yields a simple proof of the shellability of interval
orders which we give below.  Recall that a bounded poset is a poset that has
a minimum element
$\hat 0$ and a maximum element $\hat 1$. If $P$ is a bounded poset then
$\bar P$ denotes the induced subposet $P\setminus \{\hat 0, \hat 1\}$.  The
length of bounded poset
$P$ is the length of the longest chain from $\hat 0$ to $\hat 1$.   For any
$a
\le b$ in $P$, the open interval
$\{x \in P \mid a < x < b\}$ is denoted by $(a,b)$ and the closed interval
$\{x \in P \mid a \le x \le b\}$ is denoted by $[a,b]$.  The order
complex of $P$ is the simplicial complex of chains of $P$
and is denoted by $\Delta(P)$.

\proclaim{Proposition 10} (Billera and Meyers)  Every interval order is
shellable. 
\endproclaim

\demo{Proof} Let $P$ be an interval order.  By the well-known
characterization of interval orders, $P$ does not contain $Q$ (the poset
with two disjoint 2 element chains) as an induced subposet.  We may assume
without loss of generality that
$P$ is bounded and that $P$ has more than one atom.  

 The fact
that
$P$ does not contain $Q$ enables us to  choose an atom
$a$  such that each of the  covers of $a$ is greater than some other
atom. Since
$\lk_{\Delta(\bar P)}(a) = \Delta((a,\hat 1))$,
this implies that no facet of
$\lk_{\Delta(\bar P)}(a)
$ is a facet of $\Delta(\bar P\setminus \{a\})$.  Also the interval
$(a,\hat 1)$  and the induced subposet $\bar P\setminus \{a\}$ both
inherit the property of not containing $Q$ as an induced subposet. 
Hence by induction they are shellable.  We conclude that
$\Delta(\bar P)$ and hence
$P$ are shellable by Lemma 6. \hfill\qed 
\enddemo

\demo{Remark} Bj\"orner \cite{B} has independently used the same idea
to prove more generally that all interval orders are vertex decomposable.
\enddemo

Another   proof that interval orders are shellable can be
obtained using the technique of lexicographical shellability, cf. \cite
{BW2}.

\proclaim{Theorem 11} Every bounded interval order is CL-shellable.
\endproclaim

\demo{Proof} Let $P$ be a bounded interval order. 
Partially order the atoms of $P$ by letting $a \prec b$ if $a$ has a
cover that is not greater than $b$.  Antisymmetry and transitivity
follow readily from the forbidden induced subposet characterization of
interval order.  It is straight forward to verify that any linear
extension of $\preceq$ is a recursive atom ordering of $P$ by
induction on $|P|$.  Therefore $P$ is CL-shellable. \hfill\qed
\enddemo

For any bounded poset $P$ let $\beta_i(P)$ be the ith reduced Betti number of
$\Delta(\bar P)$.  It is easy to see that in
the partial order on atoms given in the proof of Theorem~11 there is a
unique minimum atom. We shall refer to this atom as the {\it
smallest} atom.

\proclaim {Corollary 12} Let $P$ be bounded  interval order of length $\ge
2$, let
$A$ be its set of atoms and let $a_0$ be its smallest atom.  Then for $i
\ge 0$,
$$\beta_i(P) = \sum_{a \in A\setminus \{a_0\}} \beta_{i-1}([a,\hat
1]).$$
\endproclaim

\demo{Proof}   By \cite{BW2, Theorem~5.9}, $\beta_i(P)$ is the number
of falling maximal chains of length $i+2$ with respect to the
CL-labeling induced by the recursive atom ordering given in Theorem~11.
So we need to describe these falling chains.  Each falling chain
from $\hat 0$ to $\hat 1$ of length $i+2$  is of the form $\{\hat
0\}
\cup c$ where
$c$ is a falling chain of length $i+1$ from $a$ to $\hat 1$ for some
atom
$a$.  We need to determine which atoms $a$ and falling chains $c$
from $a$ to $\hat 1$ are such that $\{\hat 0\} \cup c$ is falling.
The proof of \cite{BW1, Theorem 3.2} produces a CL-labeling from a
recursive atom ordering (although it's done in the pure case in
\cite{BW1}, it easily carries over to general case, cf. \cite{BW2}).   A
maximal chain has a descent on the subchain $\hat 0 \to a \to b$
if and only if
$b$ is greater than an atom that precedes
$a$ in the recursive atom ordering.  This happens for every
maximal chain through
 $a
\ne a_0$ and for no maximal chain through $a_0$.  Hence the
maximal chains of the form $\{a\} \cup c$, where $a \ne a_0$ and $c$ is
a falling chain of
$[a,\hat 1]$, are the falling chains of $P$.\hfill\qed
\enddemo

The problem of studying obstructions could conceivably be made easier by
considering special classes of simplicial complexes that are closed
under taking induced subcomplexes.  A natural class, suggested by Bj\"orner
\cite{B}, is the class of flag complexes.  One might ask whether the pair of
disjoint edges is the only obstruction for flag complexes.  It turns out
that this is not the case.  The obstruction
$M_7$ given in the proof of Theorem~5 is a flag complex.  However
obstructions $M_5$ and
$M_6$ are not flag complexes.  Also the obstructions given in
the proof of Proposition 1 are not flag complexes. So  there might still be
something interesting to say about flag complex obstructions.  

\medspace

\heading Acknowledgment \endheading

I am    grateful to Richard Stanley 
for  stimulating discussions which initiated this study of 
obstructions.

\enddocument